\documentclass[conference]{IEEEtran}
\IEEEoverridecommandlockouts
\usepackage{times}
\pdfoutput=1
\let\chapter\section
\renewcommand\dbltextfloatsep{3pt}
\renewcommand\dblfloatsep{3pt}
\renewcommand\textfloatsep{3pt}
\usepackage{times}

\usepackage[numbers]{natbib}
\usepackage{multicol}
\usepackage[bookmarks=true]{hyperref}
\usepackage{amsfonts}       

\usepackage{algpseudocode,algorithm,algorithmicx}
\usepackage{amsmath}
\usepackage{subfigure, graphicx}
\newtheorem{remark}{Remark}

\newtheorem{Theorem}{Theorem}

\allowdisplaybreaks

\begin{document}
\bstctlcite{IEEEexample:BSTcontrol}
\title{A Separation-based Approach to Data-based Control for Large-Scale Partially Observed Systems}

\author
{
Dan Yu$^{1}$, Mohammadhussein Rafieisakhaei$^{2}$ and Suman Chakravorty$^{1}$
\thanks{*This material is based upon work partially supported by NSF under Contract Nos. CNS-1646449 and Science \& Technology Center Grant CCF-0939370, the U.S. Army Research Office under Contract No. W911NF-15-1-0279, and NPRP grant NPRP 8-1531-2-651 from the Qatar National Research Fund, a member of Qatar Foundation, AFOSR contract Dynamic Data Driven Application Systems (DDDAS) contract FA9550-17-1-0068 and NSF NRI project ECCS-1637889.}
\thanks{$^{1}$D. Yu and S. Chakravorty are with the Department of Aerospace Engineering, and $^{2}$M. Rafieisakhaei is with the Department of Electrical and Computer Engineering, Texas A\&M University, College Station, Texas, 77840 USA.
        \{\tt\small yudan198811@hotmail.com, mrafieis, schakrav@tamu.edu\}}%
}

\maketitle
\begin{abstract}
This paper studies the partially observed stochastic optimal control problem for systems with state dynamics governed by partial differential equations (PDEs) that leads to an extremely large problem. First, an open-loop deterministic trajectory optimization problem is solved using a black-box simulation model of the dynamical system. Next, a Linear Quadratic Gaussian (LQG) controller is designed for the nominal trajectory-dependent linearized system which is identified using input-output experimental data consisting of the impulse responses of the optimized nominal system. A computational nonlinear heat example is used to illustrate the performance of the proposed approach. 
\end{abstract}

\section{Introduction}\label{Section 1}
Stochastic optimal control problems, also known as Markov decision problems (MDPs), have found numerous applications in the Sciences and Engineering. In general, the goal is  to control a stochastic system subject to transition uncertainty in the state dynamics so as to minimize the expected running cost of the system. The MDPs are termed Partially Observed (POMDP) if there is sensing uncertainty in the state of the system in addition to the transition uncertainty. In this paper, we consider the stochastic control of partially observed nonlinear dynamical systems that are governed by partial differential equations (PDE). In particular, we propose a novel data based approach to the solution of very large POMDPs wherein the underlying state space is obtained from the discretization of a PDE: problems whose solution has never been hitherto attempted using approximate MDP solution techniques.  

It is well known that the global optimal solution for MDPs can be found by solving the Hamilton-Jacobi-Bellman (HJB) equation \cite{dp_bertsekas}. The solution techniques can be further divided into model based and model free techniques, according as whether the solution methodology uses an analytical model of the system or it uses a black box simulation model, or actual experiments.

In model based techniques, many methods \cite{dp_num} rely on a discretization of the underlying state and action space, and hence, run into the "curse of dimensionality (COD)", the fact that the computational complexity grows exponentially with the dimension of the state space of the problem.  The most computationally efficient among these techniques are trajectory-based methods, first described in \cite{bryson}. These methods expand the nonlinear system equations about a deterministic nominal trajectory, and perform a localized version of policy iteration to iteratively optimize the trajectory. For example, the differential dynamic programming (DDP) \cite{ddp, sddp} linearizes the dynamics and the cost-to-go function around a given nominal trajectory, and designs a local feedback controller using DP. The iterative Linear Quadratic Gaussian  (ILQG) \cite{ilqg1, ilqg2}, which is closely related to DDP,  considers the first order expansion of the dynamics (in DDP, a second order expansion is considered), and designs the feedback controller using Riccati-like equations, and is shown to be computationally more efficient. In both approaches, the control policy is executed to compute a new nominal trajectory, and the procedure is repeated until convergence. 

In the model free solution of MDPs, the most popular approaches are the adaptive dynamic programming  (ADP) \cite{adp,adp_3} and  reinforcement learning (RL) paradigms \cite{rl_1, rl_3}. They are essentially the same in spirit, and seek to improve the control policy for a given black box system by repeated interactions with the environment, while observing the system's responses. The repeated interactions, or learning trials, allow these algorithms to construct a solution to the DP equation, in terms of the cost-to-go function,  in an online and recursive fashion. Another variant  of RL techniques is the so-called Q-learning method, and the basic idea in Q-learning is to estimate a real-valued function $Q(x, a)$ of states and actions instead of the cost-to-go function $V(x)$. For continuous state and control space problems, the cost-to-go functions and the Q-functions are usually represented in a functionally parameterized form, for instance,  in the linearly parametrized form $Q(x,a) = \theta' \phi(x,a)$,  where $\theta$ is the unknown parameter vector, and $\phi$ is  a pre-defined basis function, $(\cdot)'$ denotes the transpose of $(\cdot)$. Multi-layer neural networks may also be used as nonlinearly parameterized approximators instead of the linear architecture above. The ultimate goal of these techniques is the estimation/ learning of the parameters $\theta$ from learning trials/ repeated simulations of the underlying system. However, the size of the parameter $\theta$ grows exponentially in the size of the state space of the problem without a compact parametrization of the cost-to-go or Q function in terms of the a priori chosen basis functions for the approximation, and hence, these techniques are typically subject to the curse of dimensionality. Albeit a compact parametrization may exist, a priori, it is usually never known.  

In the past several years, techniques based on the differential dynamic programming/ ILQG approach \cite{ddp,sddp, RLHD4,RLHD5}, such as the RL techniques \cite{RLHD1, RLHD2} have shown the potential for RL algorithms to scale to higher dimensional continuous state and control space problems, in particular, high dimensional robotic task planning and learning problems. These methods are a localized version of the policy search \cite{baxter, sutton3,marbach, Survey_PS} technique that seek to directly optimize the feedback policy via a compact parameterization. For continuous state and control space problems, the method of choice is to wrap an LQR feedback policy around a nominal trajectory and then perform a recursive optimization of the feedback law, along with the underlying trajectory, via repeated simulations/ iterations. However, the parametrization can still be very large for partially observed problems (at least $O(d^2))$ where $d$ is the dimension of the state space) or large motion planning problems such as systems governed by partial differential equations wherein the (discretized) state is very high dimensional (thousands/ millions of states) which are typically partially observed thereby compounding the problem. Furthermore, there may be convergence problems with these techniques that can lead to the so-called ``policy chatter" phenomenon \cite{RLHD1}. 

\textit{Fundamentally, rather than solve the derived ``Dynamic Programming" problem as in the majority of the approaches above  that requires the optimization of the feedback law, our approach is to directly solve the original stochastic optimization problem in a ``separated open loop -closed loop" fashion wherein:  1) we solve an open loop deterministic optimization problem to obtain an optimal nominal trajectory in a model free fashion, and then 2) we design a closed loop controller for the resulting linearized time-varying system around the optimal nominal trajectory, again in a model free fashion. Nonetheless, the above ``divide and conquer" strategy can be shown to be near optimal. }  

The primary contributions of the proposed approach are as follows:

1) We specify a detailed set of experiments to accomplish the closed loop controller design for any unknown nonlinear system, no mater how high dimensional. This series of experiments consists of a sequence of input perturbations to collect the impulse responses of the system, first to find an optimized nominal trajectory, and then to recover the LTV system corresponding to the perturbations of the nominal system in order to design an LQG controller for the LTV system.

2) In general, for large scale systems with partially observed states, the system identification algorithm such as time-varying ERA \cite{tv_era} automatically constructs reduced order model (ROM) of the LTV system, and hence, results in a reduced order estimator and controller. Therefore, even for large scale systems such as partially observed systems with the state dynamics governed by PDEs, the computation of the feedback controller is nevertheless computationally tractable, for instance, in the partially observed nonlinear heat control problem considered in this paper, the complexity is reduced by $O(10^5)$ when compared to DDP based RL techniques.

3) We provide a unification of traditional linear and nonlinear optimal control techniques with ADP and RL techniques in the context of Stochastic Dynamic Programming problems.

The rest of the paper is organized as follows. In Section \ref{Section 2}, the basic problem formulation is outlined. In Section \ref{Section 3}, we propose a separation based stochastic optimal control algorithm, with discussions of implementation problems. In Section \ref{Section 4}, we  test the proposed approach using a one-dimensional nonlinear heat problem.

\section{Problem Setup}\label{Section 2}
Consider a discrete time nonlinear dynamical system:
\begin{align} \label{original system}
x_{k + 1} &= f(x_k, u_k, w_k), \nonumber \\
y_k &= h(x_k , v_k),
\end{align}
where $x_k \in \mathbb{R}^{n_x}$,  $y_k \in  \mathbb{R}^{n_y}$, $u_k \in \mathbb{R}^{n_u}$ are the state vector, the measurement vector and the control vector at time $k$ respectively. The system function $f(\cdot)$ and measurement function $h(\cdot)$ are nonlinear. The process noise $w_k$ and measurement noise $v_k$ are assumed as zero-mean, uncorrelated Gaussian white noise, with covariance $W$ and $V$ respectively. In considering PDEs, the dynamics above are the discretized version of the equations (using Finite Difference (FD) or Finite Element (FE) schemes). Typically, the discretization leads to a very large state space problem consisting of at least hundreds of states and typically millions of states for larger problems.

The belief $b(x_k)$ is defined as the distribution of the state $x_k$ given all past control inputs and sensor measurements, and is denoted by $b_k$. In this paper, we represent beliefs by Gaussian distributions, and denote the belief $b_k = (\mu_k, \Sigma_k)$, where $\mu_k$  and $\Sigma_k$ are the mean and covariance of the Gaussian belief state. Denote the belief dynamics 
\begin{align}
b_{k + 1} = \tau (b_k, u_k, y_{k + 1}),
\end{align}
and assume that $b_0$ is known. Note that if the belief is Gaussian, the belief state is $O(n_x^2)$, which for a PDE is extremely large due to the fact that $n_x$ is very large. 

In this paper, we consider the following stochastic optimal control problem. 

\textbf{Stochastic Control Problem:} For the system with unknown nonlinear dynamics, i.e., $f(\cdot)$ and $h(\cdot)$ are unknown, the optimal control problem is to find the control policies $\pi = \{\pi_0, \pi_1, \cdots, \pi_{N -1} \}$ in a finite time horizon $[0, N]$, where $\pi_k$ is the control policy at time $k$, i.e., $u_k = \pi_k(b_k)$, such that for a given initial belief state $b_0$,  the cost function
\begin{align}\label{cost_sto_orig}
J_{\pi} = E(\sum_{k = 0}^{N - 1} c_k(b_k, u_k) + c_N(b_N)), 
\end{align}
is minimized, where $\{ c_k (\cdot,\cdot) \}_{k = 0}^{N - 1}$ denotes the  immediate cost function, and $c_N(\cdot)$ denotes the terminal cost. The expectation is taken over all randomness. 

\section{Separation based Feedback Control Design}\label{Section 3}

The stochastic control problem is solved in a separated open loop- closed loop (SOC) fashion, i.e., first, we solve a noiseless open-loop optimization problem to find a nominal optimal trajectory and then we design a linearized closed-loop controller around the nominal trajectory, such that, with existence of stochastic perturbations, the state stays close to the optimal open-loop trajectory. The three-step framework to solve the stochastic feedback control problem may be summarized as follows. 
\begin{itemize}
\item Solve the open loop optimization problem using a general nonlinear programming (NLP) solver with a black box simulation model of the dynamics, where the belief dynamics is updated using an Ensemble Kalman Filter (EnKF) \cite{enkf1}.
\item Linearize the system around the nominal open loop optimal  belief  trajectory, and identify the linearized time-varying system from input-output experiment data using a suitable system identification algorithm such as the time-varying eigensystem realization algorithm (ERA) \cite{tv_era}.
\item Design an LQG controller which results in an optimal linear control policy around the nominal trajectory.
\end{itemize}
In the following section, first, we present the ``Separation" theorem and then discuss each of the above steps.

\subsection{A Separation Result}\label{sep_res}
\textit{Nominal trajectories:} Denote $ \{{\bar{u}}_k\}_{k=0}^{N-1},$ $\{{\bar{\mu}}_k\}_{k=0}^{N}$ as the nominal control and state trajectories of the system, respectively, $\{{\bar{y}}_k\}_{k=0}^{N}$ as the corresponding observations and $\{{\bar{b}}_k\}_{k=0}^{N} $ as the belief trajectories, where:
\begin{align}\label{eq:Nominal trajectories}
&{\bar{\mu}}_{k+1}={f}({\bar{\mu}}_{k},{\bar{u}}_k, 0),
{\bar{y}}_{k}={h}({\bar{\mu}}_{k}, 0), \nonumber \\
& {\bar{u}}_{k}={\pi}_{k}({\bar{b}}_k),
{\bar{b}}_{k + 1} = {\tau} ({\bar{b}}_k, {\bar{u}}_k, {\bar{y}}_{k + 1}),
\end{align}
with the initial conditions of $ {\bar{b}}_0 = {b}_0 $, and ${\bar{\mu}}_{0}=\mathbb{E}[{b}_0] $.


\textit{Nominal cost function:} The nominal cost and its first order expansion are given by (please see \cite{Separation} for details):
\begin{align}
\bar{J}:=&\sum_{k=0}^{N-1}c_k({\bar{b}}_k, {\bar{u}}_k)+c_N({\bar{b}}_N), 
\\J\approx& \bar{J} + \underbrace{\sum_{k=0}^{N-1} ({C}^{{b}}_k({b}_k-{\bar{b}}_k)+ {C}^{{u}}_k({u}_k-{\bar{u}}_k))+ {C}^{{b}}_K({b}_N-{\bar{b}}_N)}_{\delta J} 
\end{align}

\begin{Theorem}[Cost Function Linearization Error]\label{theorem:First Order Cost Function Error}  The expected first-order linearization error of the cost function is zero, $\mathbb{E}(\delta J) = 0$. 
\end{Theorem}
A typical stochastic trajectory optimization consists of optimizing the nominal trajectory along with an associated linearized feedback controller \cite{ddp, RLHD1, RLHD2}. Theorem \ref{theorem:First Order Cost Function Error} shows that the  first order approximation of the stochastic cost function is dominated by the nominal cost and depends only on the nominal trajectories of the system, which is independent of the linear feedback controller designed to track the optimized nominal system. Therefore, the design of the optimal feedback gain can be separated from the design of the optimal nominal trajectory of the system. As a result, the stochastic optimal control problem can be divided into two separate problems: the first is a deterministic problem to design the open-loop optimal control sequence, and hence, the optimal nominal trajectory of the system. The second problem is the design of an optimal linear feedback law to track the nominal trajectory of the system. Note that in the case of a belief space problem, the nominal trajectory is the optimal belief state trajectory unlike typical trajectory optimization based RL methods designed for fully observed problems such as \cite{RLHD1, RLHD2}.


%

\subsection{Open Loop Trajectory Optimization in Belief Space}
Consider the following open loop belief state optimization problem with given initial belief state $b_0$:
\begin{align}\label{cost_open}
\{ u_k^*\}_{k = 0}^{N - 1} = &\null  \arg \min_{\{u_k\}} \bar{J} (\{b_k\}_{k=0}^N, \{u_k\}_{k=0}^{N-1}), \nonumber \\
&\null b_{k + 1} = \tau(b_k, u_k, \bar{y}_{k + 1}),
\end{align}
where the nominal observations $\bar{y}_k$ are generated as follows:
\begin{align}
x_{k + 1} = f(x_k, u_k, 0), \bar{y}_k = h(x_k, 0),
\end{align}
with $x_0 = \mu_0$. Note that given the nominal observations $\bar{y}_k$, the belief evolution is deterministic and hence, the above is a deterministic optimization problem (this was first posed in the reference \cite{platt10} in the context of an Extended Kalman Filter (EKF) based belief propagation scheme).
\subsubsection{Belief Propagation using EnKF} Given an initial belief state $b_0 = (\mu_0, \Sigma_0)$, and a control sequence $\{ u_k \}$, the EnKF algorithm can be used to propagate the belief state using only a simulator of the state dynamics. This is necessary since we typically do not have access to even an (approximate) belief state simulator for large scale systems such as PDEs, and hence, need to construct one from a state space simulator. Note that the state space dimension is at least in the hundreds for such systems, and thus, a particle filter would suffer from particle depletion, and hence, cannot be used.

Denote the EnKF algorithm as
\begin{align}
\{ b_k \}_{ k = 0}^{N} = \textbf{EnKF} (b_0, \{ u_k \}_{k = 0}^{N - 1}).
\end{align}
The details of EnKF can be found in \cite{enkf2}, and is briefly summarized in Appendix \ref{app_enkf}: in short, it is a particle filter that is free from the particle depletion problem and is typically used in the filtering of large scale systems such as those governed by PDEs, for instance, meteorological phenomenon.

\subsubsection{Open loop optimization approach} The open loop optimization problem is solved using the gradient descent approach \cite{gradient, sim_opt} utilizing an EnKF. Denote  the initial guess of the control sequence as $U^{(0)} = \{u_k^{(0)} \}_{k = 0}^{N - 1}$, and  the corresponding belief state $\mathcal{B}^{(0)} = \{b_k^{(0)} \}_{k=0}^N = \textbf{EnKF}(b_0, U^{(0)})$.

The control policy is updated iteratively via
\begin{align}
U^{(n + 1)} = U^{(n)} - \alpha \nabla_U \bar{J}(\mathcal{B}^{(n)}, U^{(n)}),
\end{align}
until a convergence criterion is met, where $U^{(n)} = \{u_k^{(n)} \}_{k = 0}^{N - 1}$ denotes the control sequence in the $n^{th}$ iteration, $\mathcal{B}^{(n)} = \{b_k^{(n)}\}_{k=0}^N$ denotes the corresponding belief state, and $\alpha$ is the step size parameter.

The gradient vector is defined as:
\begin{align}\label{gradv}
\nabla_U \bar{J}(\mathcal{B}^{(n)}, U^{(n)}) = \begin{pmatrix} \frac{\partial \bar{J}}{\partial u_0} & \frac{\partial \bar{J}}{\partial u_1} & \cdots & \frac{\partial \bar{J}}{\partial u_{N - 1}} \end{pmatrix}|_{\mathcal{B}^{(n)}, U^{(n)}},
\end{align}
and without knowing the explicit form of the cost function, each partial derivative with respect to the $i^{th}$ control variable $u_i$ is calculated as follows:
\begin{align}\label{pd}
\frac{\partial \bar{J}}{\partial u_i}|_{\mathcal{B}^{(n)}, U^{(n)}} = \frac{1}{h}(\bar{J}(\mathcal{B}_i^{(n)}, u_0^{(n)},  \cdots, u_i^{(n)} + h, \cdots, u_{N - 1}^{(n)}) \nonumber\\
- \bar{J}(\mathcal{B}^{(n)}, u_0^{(n)},  \cdots, u_i^{(n)}, \cdots, u_{N - 1}^{(n)})),
\end{align}
where $h$ is a small constant perturbation and $\mathcal{B}_i^{(n)}$ denotes the belief state corresponding to the control input
$\{ u_0^{(n)},  \cdots, u_i^{(n)} + h, \cdots, u_{N - 1}^{(n)} \}$, $i = 0, \cdots, N-1$.

The open loop optimization approach is summarized in Algorithm \ref{algo_grad}.

\begin{algorithm}[!bt]
\caption{Open Loop Optimization Algorithm}
\label{algo_grad}
 \begin{algorithmic}[1]
 \Require{Start belief $b_0$, cost function $\bar{J}(.)$, initial guess $U^{(0)} = \{u_k^{(0)} \}_{k = 0}^{N - 1}$,  gradient descent design parameters  $\alpha, h, \epsilon$.}
 \Ensure{Optimal control sequence $\{ \bar{u}_k \}_{k = 0}^{N - 1}$, belief nominal trajectory $\{\bar{b}_k \}_{k  = 0}^N$}
 \State $n = 0$, set $\nabla_U \bar{J}(\mathcal{B}^{(0)}, U^{(0)}) = \epsilon $.
 \While{$\nabla_U \bar{J}(\mathcal{B}^{(n)}, U^{(n)}) \geq \epsilon $}
 \State Compute the belief $\mathcal{B}^{(n)} = $ \textbf{EnKF}($b_0, U^{(n)}$).
 \State Evaluate the cost function $\bar{J}(\mathcal{B}^{(n)}, U^{(n)})$.
 \State Perturb each control variable $u_i^{(n)}$ by $h$ and compute the belief $\mathcal{B}_i^{(n)}$, $i = 0, \cdots, N - 1$, calculate the gradient vector $\nabla_U \bar{J}(\mathcal{B}^{(n)}, U^{(n)})$.
  \State Update the control policy
 $U^{(n + 1)} = U^{(n)} - \alpha \nabla_U \bar{J}(\mathcal{B}^{(n)}, U^{(n)}).$
 \State $n = n + 1$.
 \EndWhile
 \State $\{\bar{u}_k \}_{k = 0}^{N - 1} =  U^{(n)}$, $\{\bar{b}_k \}_{k  = 0}^N$ = \textbf{EnKF}($b_0, U^{(n)}$).
\end{algorithmic}
\end{algorithm}

\begin{remark}
The open loop optimization problem is solved using a black box simulation model of the underlying dynamics, with a sequence of input perturbation learning trials. Higher order approaches other than gradient descent are possible \cite{sim_opt}, however, for a general system, the gradient descent approach is easy to implement, and is amenable to very large scale parallelization.
\end{remark}

\subsection{Linear Time-Varying System Identification}
Denote the optimal open-loop control  as $\{\bar{u}_k \}_{k = 0}^{N - 1}$, and the corresponding nominal belief state as $\{ \bar{\mu}_k, \bar{\Sigma}_k\}_{k = 0}^N$.  We linearize the system (\ref{original system}) around the nominal trajectory (the mean $\{\bar{\mu}_k\}$), assuming that the control and disturbance enter through the same channels and the noise is purely additive (these assumptions are only for simplicity and can be relaxed easily):
\begin{align}\label{perturbation system}
& \delta x_{k+ 1} = A_k \delta x_k + B_k (\delta u_k + w_k), \nonumber \\
& \delta y_k = C_k \delta x_k +  v_k,
\end{align}
where $\delta x_k = x_k  - \bar{\mu}_k$ describes the state deviations from the nominal mean trajectory, $\delta u_k = u_k - \bar{u}_k$ describes the control deviations,  $\delta y_k = y_k - h(\bar{\mu}_k, 0)$ describes the measurement deviations, and
\begin{align}
& A_k = \frac{\partial f(x, u, w)}{\partial x}|_{\bar{\mu}_k, \bar{u}_k, 0}, B_k = \frac{\partial f(x, u, w)}{\partial u}|_{\bar{\mu}_k, \bar{u}_k, 0},   \nonumber\\
& C_k = \frac{\partial h(x, v)}{\partial x}|_{\bar{\mu}_k, 0}.
\end{align}

Consider system (\ref{perturbation system}) with zero noise and  $\delta x_0 = 0$,  the input-output relationship is given by:
\begin{align}\label{input_output}
\delta y_k = \sum_{j = 0}^{k - 1} h_{k, j} \delta u_j,
\end{align}
where $h_{k, j}$ is defined as the generalized Markov parameters, and
\begin{align}
h_{k, j} \begin{cases} = C_k A_{k - 1} A_{k - 2} \cdots A_{j+1} B_j, & \text{ if } j < k - 1, \\
= C_k B_{k -1},  & \text{ if } j = k - 1, \\
 = 0, & \text{ if } j > k - 1.
 \end{cases}
\end{align}

\subsubsection{\textbf{Partial Realization Problem \cite{antoulas, skelton1}}} Given a finite sequence of Markov parameters $h_{k, j} \in \Re^{n_y \times n_u}, k = 1, 2, \cdots, s,  j = 0, 1, \cdots, k$, the partial realization problem consists of finding a positive integer $n_r$ and LTV system $(\hat{A}_k, \hat{B}_k, \hat{C}_k)$, where $\hat{A}_k \in \Re^{n_r \times n_r}, \hat{B}_k \in \Re^{n_r \times n_u}, \hat{C}_k \in \Re^{n_y \times n_r}$, such that the identified generalized Markov parameters $\hat{h}_{k, j} \equiv  \hat{C}_k \hat{A}_{k - 1} \hat{A}_{k - 2} \cdots \hat{A}_{j+1} \hat{B}_j= h_{k, j}$. Then $(\hat{A}_k, \hat{B}_k, \hat{C}_k)$ is called a partial realization of the sequence $h_{k, j}$.

We solve the partial realization problem using the time-varying ERA. Time-varying ERA starts by estimating the generalized Markov parameters using input-output experiments, constructs a generalized Hankel matrix, and solves the singular value decomposition (SVD) problem of the constructed Hankel matrix. The details of the time-varying ERA can be found in \cite{tv_era}, and is briefly summarized here.

Define the generalized Hankel matrix as:
\begin{align}\label{hankel}
\underbrace{H_k^{(p, q)}}_{pn_y \times qn_u} = \begin{pmatrix}  h_{k, k-1} & h_{k, k-2} & \cdots & h_{k, k - q} \\
h_{k + 1, k - 1} & h_{k + 1, k -2} & \cdots & h_{k + 1, k-q} \\
\vdots & \vdots & \cdots & \vdots \\
h_{k+ p- 1, k -1} & h_{k + p - 1, k -2} & \cdots & h_{k + p - 1, k -q} \end{pmatrix},
\end{align}
where $p$ and $q$ are design parameters could be tuned for best performance. Denote the rank of the Hankel matrix $H_k^{(p, q)}$ is $n_r$, then $pn_y \geq n_r, qn_u \geq n_r$.

Given the generalized Markov parameters, we construct two Hankel matrices $H_{k}^{(p, q)}$ and $ H_{k + 1}^{(p, q)}$, and then solve the singular value decomposition problem:
\begin{align}\label{svd}
H_k^{(p, q)} = \underbrace{U_k \Sigma_k^{1/2}}_{O_k^{(p)}} \underbrace{ \Sigma_k^{1/2} V_k'}_{R_{k-1}^{(q)}},
\end{align}
where the rank of the Hankel matrix $H_k^{(p, q)}$ is $n_r$ and $n_r \leq n_x$. Then
$\Sigma_k \in \mathbb{R}^{n_r \times n_r}$ is the collection of all non-zero singular values,  and $U_k \in \mathbb{R}^{p n_y \times n_r}$, $V_k \in \mathbb{R}^{q n_u \times n_r}$ are the corresponding left and right singular vectors.

Similarly, $H_{k + 1}^{(p, q)} = O_{k + 1}^{(p)} R_k^{(q)}$.

Then the identified system using time-varying ERA is:
\begin{align}\label{ltv_id}
&\underbrace{\hat{A}_k}_{n_r \times n_r} = (O_{k + 1}^{(p) \downarrow})^{+} O_k^{(p)\uparrow} \nonumber \\
&\underbrace{\hat{B}_k}_{n_r \times n_u} = R_k^{(q)}(:, 1: n_u), \nonumber \\
&\underbrace{\hat{C}_k}_{n_y \times n_r} = O_k^{(p)}(1: n_y, :),
\end{align}
where $(.)^{+}$ denotes the pseudo inverse of $(.)$, $O_{k + 1}^{(p) \downarrow}$ contains the first $(p- 1) n_y$ rows of $O_{k + 1}^{(p)}$, and $O_k^{(p)\uparrow}$ contains the last $(p-1) n_y$ rows of $O_k^{(p)}$. Here, we assume that $n_r$ is constant through the time period of interest, which could also be relaxed.

\subsubsection{Identify the Generalized Markov Parameters using Input-Output Experiments} Now the problem is how to estimate the generalized Markov parameters. Consider  the input-output map for system (\ref{perturbation system}) with zero noise and  $\delta x_0 = 0$:
\begin{align}\label{input_output}
\delta y_k = \sum_{j = 0}^{k - 1} h_{k, j} \delta u_j.
\end{align}

We run $M$ simulations and in the $i^{th}$ simulation, choose input sequence $\{ \delta u_{t, (i)} \}_{t = 0}^k$, and collect the output $\delta y_{k, (i)}$. The subscript $(i)$ denotes the experiment number. Then the generalized Markov parameters $\{ h_{k, j} \}_{j = 0}^k $ could be recovered via solving the least squares problem:
\begin{align}\label{gmarkov}
&\begin{pmatrix} \delta y_{k, (1)} & \delta y_{k, (2)} \cdots & \delta y_{k, (M)} \end{pmatrix} \nonumber \\
& = \begin{pmatrix} 0 & h_{k, k-1} & h_{k, k- 2} & \cdots & h_{k, 0} \end{pmatrix}  \times \nonumber\\
& \begin{pmatrix} \delta u_{k, (1)} &
\delta u_{k, (2)} & \cdots & \delta u_{k, (M)} \\ \delta u_{k -1, (1)}& \delta u_{k -1, (2)} & \cdots & \delta u_{k - 1, (M)} \\ \vdots & \vdots & &\vdots \\ \delta u_{0, (1)} & \delta u_{0, (2)} & \cdots & \delta u_{0, (M)} \end{pmatrix},
\end{align}
where $M$ is a design parameter and is chosen such that the least squares solution is possible.

Notice that we cannot perturb the system (\ref{perturbation system}) directly. Instead, we identify the generalized Markov parameters as follows.

Run M parallel simulations with the noise-free system:
\begin{align}\label{noise_free_input}
&x_{k + 1, (i)} = f(x_{k, (i)}, \bar{u}_k + \delta u_{k, (i)}, 0), \nonumber \\
&y_{k, (i)} = h(x_{k, (i)}, 0),
\end{align}
where $i = 1, 2, \cdots, M$, and therefore,
\begin{align}\label{noise_free_output}
\delta y_{k, (i)} = y_{k, (i)} - h(\bar{\mu}_k, 0).
\end{align}
where $(\bar{u}_k, \bar{\mu}_k)$ is the open loop optimal trajectory. Then solve the same least squares problem with (\ref{gmarkov}).



The time-varying ERA used in this paper is summarized in Algorithm \ref{algo_ltv_era}.  Denote the identified deviation system
\begin{align}\label{rom}
&\delta a_{k  +1} = \hat{A}_k \delta a_k + \hat{B}_k (\delta u_k +  w_k),  \nonumber\\
& \delta y_k = \hat{C}_k \delta a_k + v_k,
\end{align}
where $\delta a_k \in \Re^{n_r}$ denotes the reduced order model (ROM) deviation states. The dimension $n_r$ of the ROM is such that $n_r << n_x$, where $n_x$ is the dimension of the state, thereby automatically providing a compact parametrization of the problem (please also see Section \ref{Section 4}).

\begin{algorithm}[!tb]
\caption{LTV System Identification}
 \label{algo_ltv_era}
  \begin{algorithmic}[1]
 \Require{Nominal Trajectory $\{u_k\}_{k = 0}^{N - 1}, \{\bar{b}_k \}_{k = 0}^N$, design parameters $M, p, q$}
 \Ensure{$\{ \hat{A}_k, \hat{B}_k, \hat{C}_k \}$}
 \State $k = 0$
 \While{$k \leq N - 1$}
 \State Identify generalized Markov parameters with input and output experimental data using
 (\ref{gmarkov}), (\ref{noise_free_input}) and (\ref{noise_free_output}).
 \State Construct the generalized Hankel matrices $H_{k}^{(p, q)}$, $H_{k + 1}^{(p, q)}$  using (\ref{hankel}).
 \State Solve the SVD problem, and construct $\{\hat{A}_k, \hat{B}_k, \hat{C}_k \}$ using (\ref{ltv_id}).
 \State $k = k + 1$.
 \EndWhile
\end{algorithmic}
\end{algorithm}

\subsection{Closed Loop Controller Design}
Given the identified deviation system (\ref{rom}), we design the closed-loop controller to  follow the optimal nominal trajectory, which is to minimize the cost function
\begin{align}
J_f = \sum_{k = 0}^{N - 1} (\delta \hat{a}_k' Q_k \delta \hat{a}_k+ \delta u_k' R_k \delta u_k ) + \delta \hat{a}_N' Q_N \delta \hat{a}_N,
\end{align}
where $\delta \hat{a}_k$ denotes the estimates of the deviation state $\delta a_k$, $Q_k, Q_N$ are positive definite, and $R_k$ is positive semi-definite. For the linear system (\ref{rom}), the ``separation principle" of linear control theory (not the Separation result of Section \ref{sep_res}) can be used \cite{dp_bertsekas}. Using this result, the design of the optimal linear stochastic controller can be separated into the decoupled design of an optimal Kalman filter and a fully observed optimal LQR controller.  The details of the design is standard \cite{dp_bertsekas} and is shown briefly in Appendix \ref{sec_lqg}.

\begin{figure*}[tb]
\centering
\includegraphics[width= 1.0\textwidth]{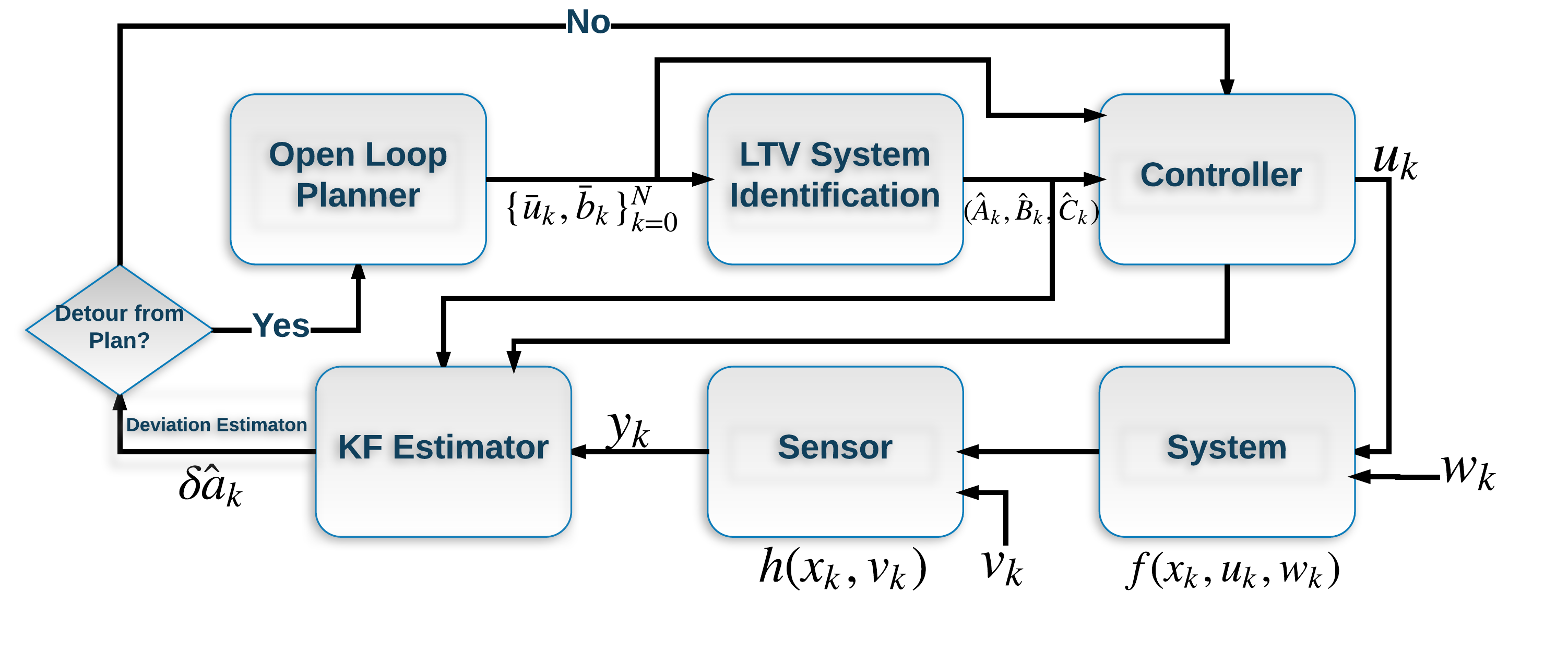}
\vspace{-0.4in}
\caption{Separation based Stochastic Feedback Control Algorithm}\label{algo_sc}
\end{figure*}

A flow chart for the Separation based Nonlinear Stochastic Control Design is shown in Fig. \ref{algo_sc}, and the algorithm is present in Algorithm \ref{algo_spc}.

\begin{algorithm}[!tb]
\caption{Separation based Stochastic Feedback Control}
\label{algo_spc}
 \begin{algorithmic}[1]
 \State Solve the deterministic open-loop optimization problem using Algorithm \ref{algo_grad}.
 \State Identify the LTV system using Algorithm \ref{algo_ltv_era}.
 \State Solve the decoupled Riccati equations (\ref{r_1}),(\ref{r_2}) using LTV system for feedback gain $\{ L_k \}_{k = 0}^{N}$.
 \State Set $k = 0$, given initial estimates $\delta \hat{a}_0 = 0, P_0$.
\While{$k \leq N - 1$}
 \State \begin{align}
& u_k = \bar{u}_k - L_k \delta \hat{a}_k, \nonumber \\
& x_{k + 1} = f(x_k, u_k, w_k), \nonumber \\
& y_{k + 1} = h(x_{k + 1}, v_{k + 1}),
\end{align}
Update $\delta \hat{a}_k$ using the Kalman Filter (\ref{k_1}), (\ref{k_2}) and (\ref{k_3}).
\State $k = k + 1$.
 \EndWhile
\end{algorithmic}
\end{algorithm}

\subsection{Discussion}
\textit{Direct Data based Controller design:}  As mentioned in data-based LQG \cite{dlqg_skelton, dlqg_2} and data-driven MPC control \cite{dd_1, dd_2}, the linear system ($\hat{A}_k, \hat{B}_k, \hat{C}_k)$ need not be identified to design the LQG controller which can be directly designed from the identified Markov parameters.


\textit{Replanning:} 
The proposed approach is theoretically valid under a small noise assumption (it is typically valid for medium noise). In practice, due to non-linearities and unknown perturbations, the actual state might deviate from the nominal trajectory during execution whence a replanning starts from the current state in a model predictive control (MPC) fashion. However, unlike in MPC, the replanning does not need to be done at every time step, only when necessary which is, in general, very infrequently.

\textit{Optimality:} The open loop law generated by the gradient descent can be guaranteed to be locally optimal under usual regularity conditions. Theorem 1 shows that, under a linear approximation, the stochastic cost is the same as the nominal cost and therefore, locally optimal as well. ILQG/ DDP based methods can also only make a claim regarding the local optimality of the nominal control law unlike global policy iteration, and therefore, the guarantees regarding optimality are the same for both.

\textit{Complexity:} The model free open loop optimization problem has complexity $O(n_u)$, where $n_u$ is the number of inputs, the LTV system identification step is again $O(n_u n_y)$, and the LQG feedback design has complexity $O(n_r^2)$, where $n_r$ is the order of the ROM from the LTV system identification step. Suppose we were to use an ILQG based design such as in \cite{RLHD1, RLHD2}, the complexity of the controller/ policy parametrization is $O(n_u n_x^2)$. Moreover, the policy evaluation step would require the estimation of a parameter of the size $O(n_x^4)$. Since $n_r<<n_x$ typically, the complexity of our separated technique is several orders of magnitude smaller (please see Section \ref{Section 4} also).

\section{Experiments}\label{Section 4}
We test the method on a one-dimensional nonlinear heat transfer problem. The heat transfer along a slab is governed by the partial differential equation:
\begin{align}
\frac{\partial T}{\partial t} = K(x, T) \frac{\partial^2 T}{\partial x^2} - \eta T + u(t),
\end{align}
where $T(x,t)$ denotes the temperature distribution at location $x$ and time $t$. The length of the slab $L = 0.6m$. $K(x, T)$ denotes the thermal diffusivity, $\eta$ denotes the convective heat transfer coefficient, and $u(t)$ denotes the external heat sources. 

The initial condition and boundary conditions are: 
\begin{align}
&T(x, 0) = 100 ^{\circ}F, \nonumber \\
&\frac{\partial T}{\partial x}|_{x = 0} = 0, T(L, t) = 150 ^{\circ}F.
\end{align}

The system is discretized using finite difference method, and there are 100 grid points which are equally spaced. We use a time step of $0.25s$. There are five point sources evenly located between $[0.1L, 0.9L]$. The sensors are placed at the same locations. Note that if we were to use an ILQG based design, the size of the state space would be 10100, and the policy evaluation step would require the solution of a 10100 x 10100 Ricatti equation.

The total simulation time is $62.5 s$. The control objective is to reach the target temperature $T_f = (150 \pm 3)^{\circ}F$ for the entire field within $t = 37.5 s$, and keep the temperature at $(150 \pm 3)^{\circ}F$ between $[37.5, 62.5]s$. 

\begin{figure*}[htbp]
\centering
\subfigure[Open Loop Nominal Belief Trajectory]{
\includegraphics[width= 0.3\textwidth]{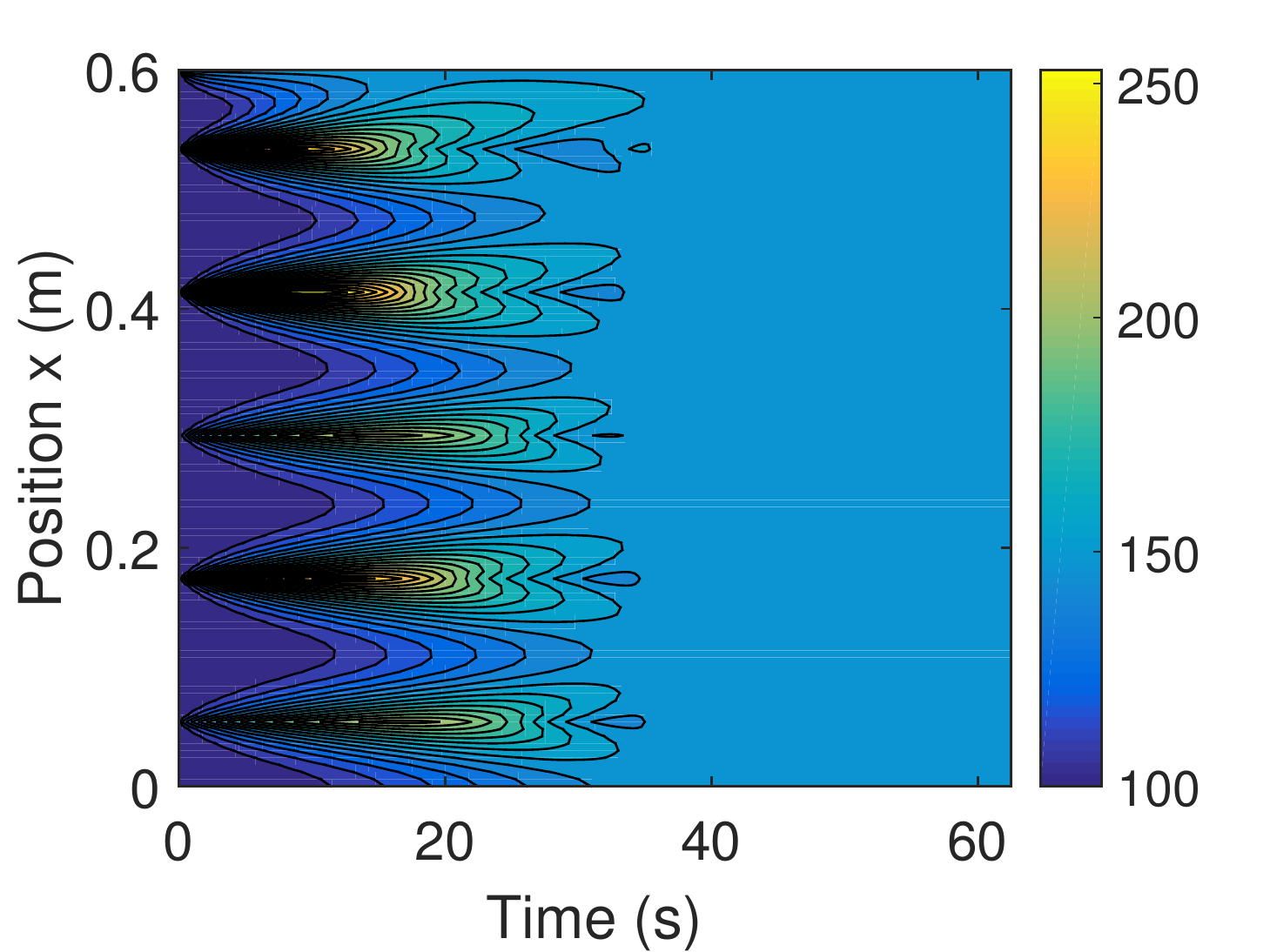}}
\subfigure[Open Loop Optimal Control]{
\includegraphics[width= 0.3\textwidth]{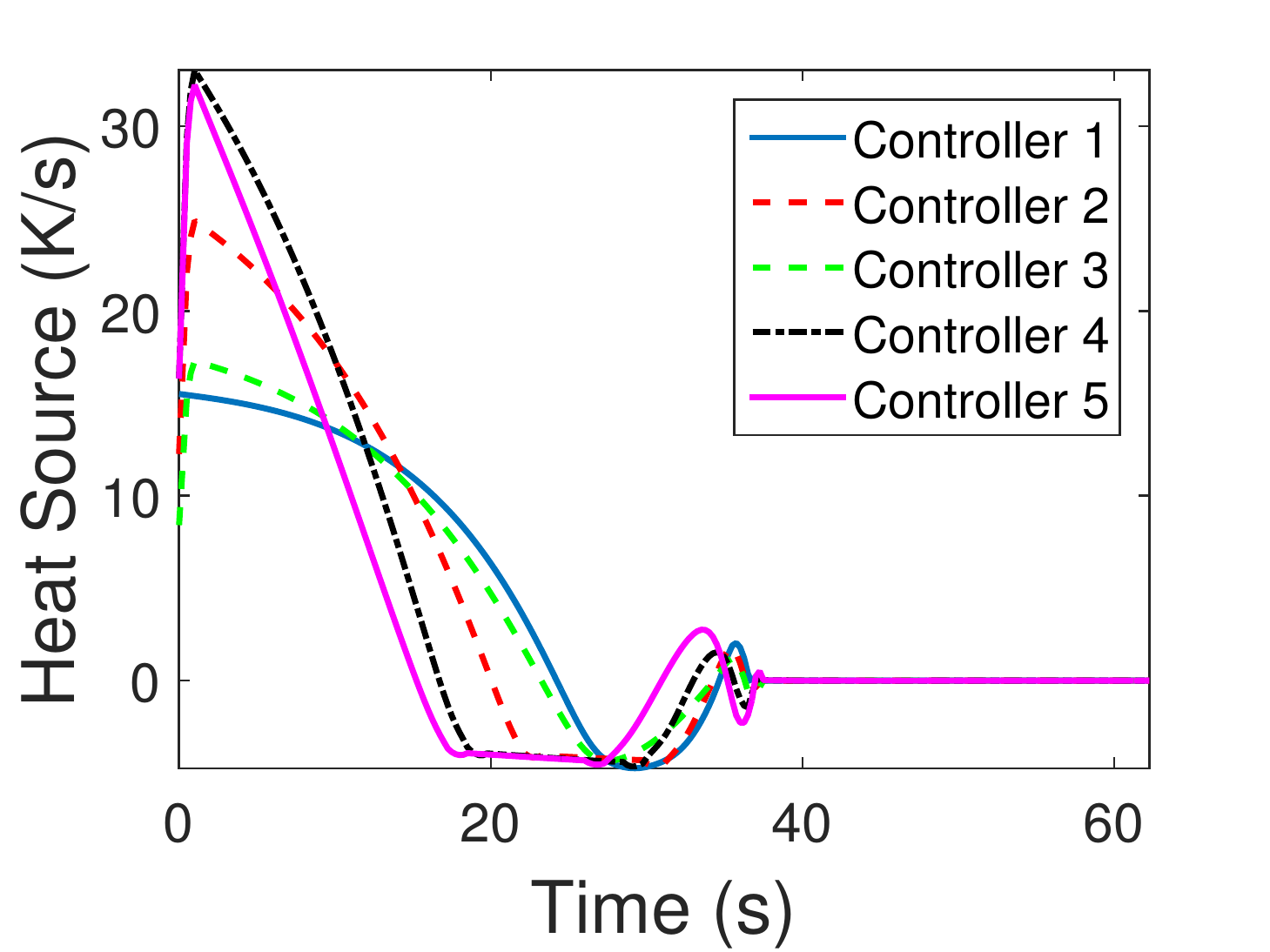}}
\subfigure[Comparison of Markov Parameters]{\includegraphics[width = 0.3\textwidth]{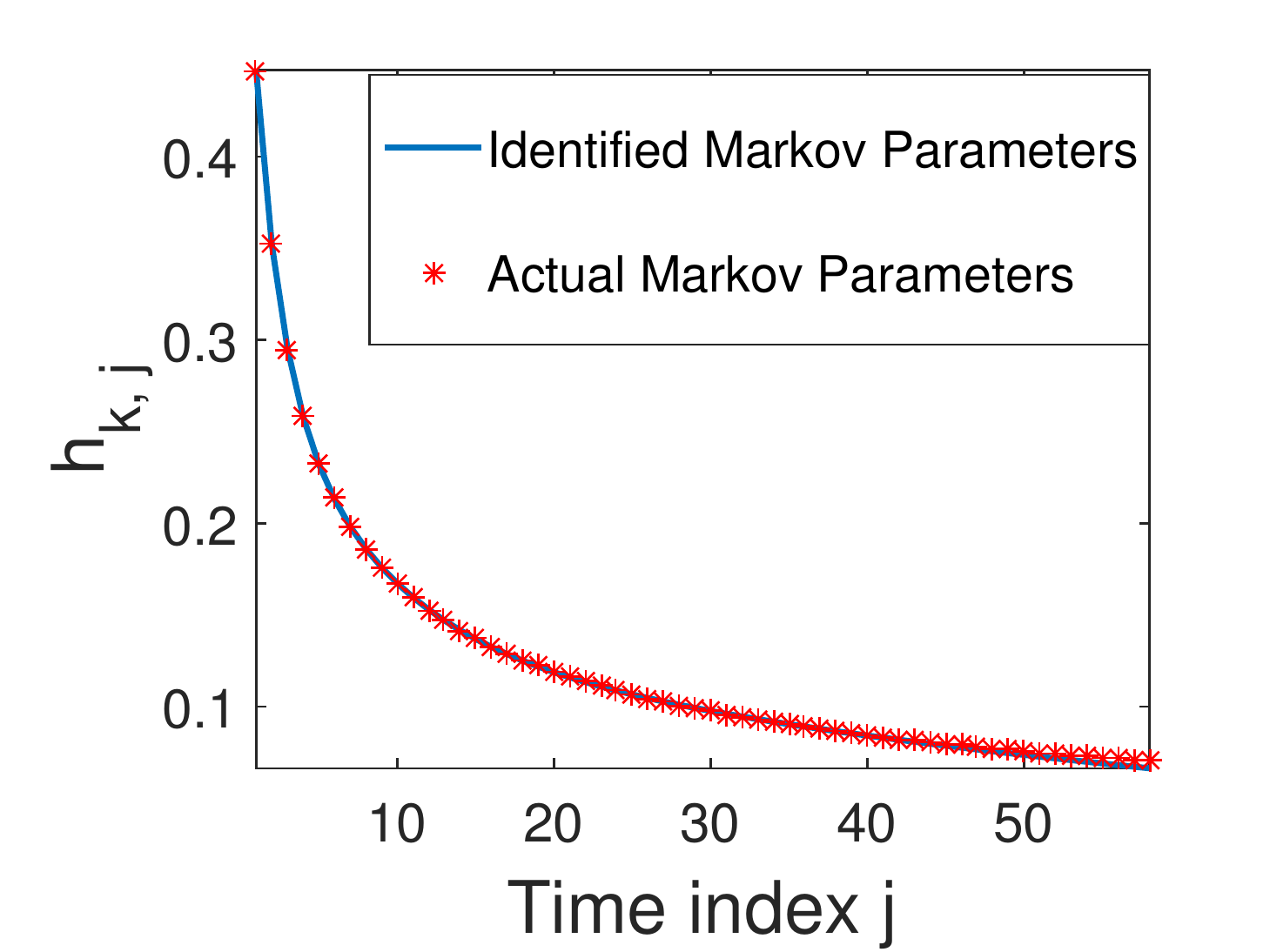}}
\subfigure[Comparison of Belief Trajectory]{
\includegraphics[width= 0.3 \textwidth]{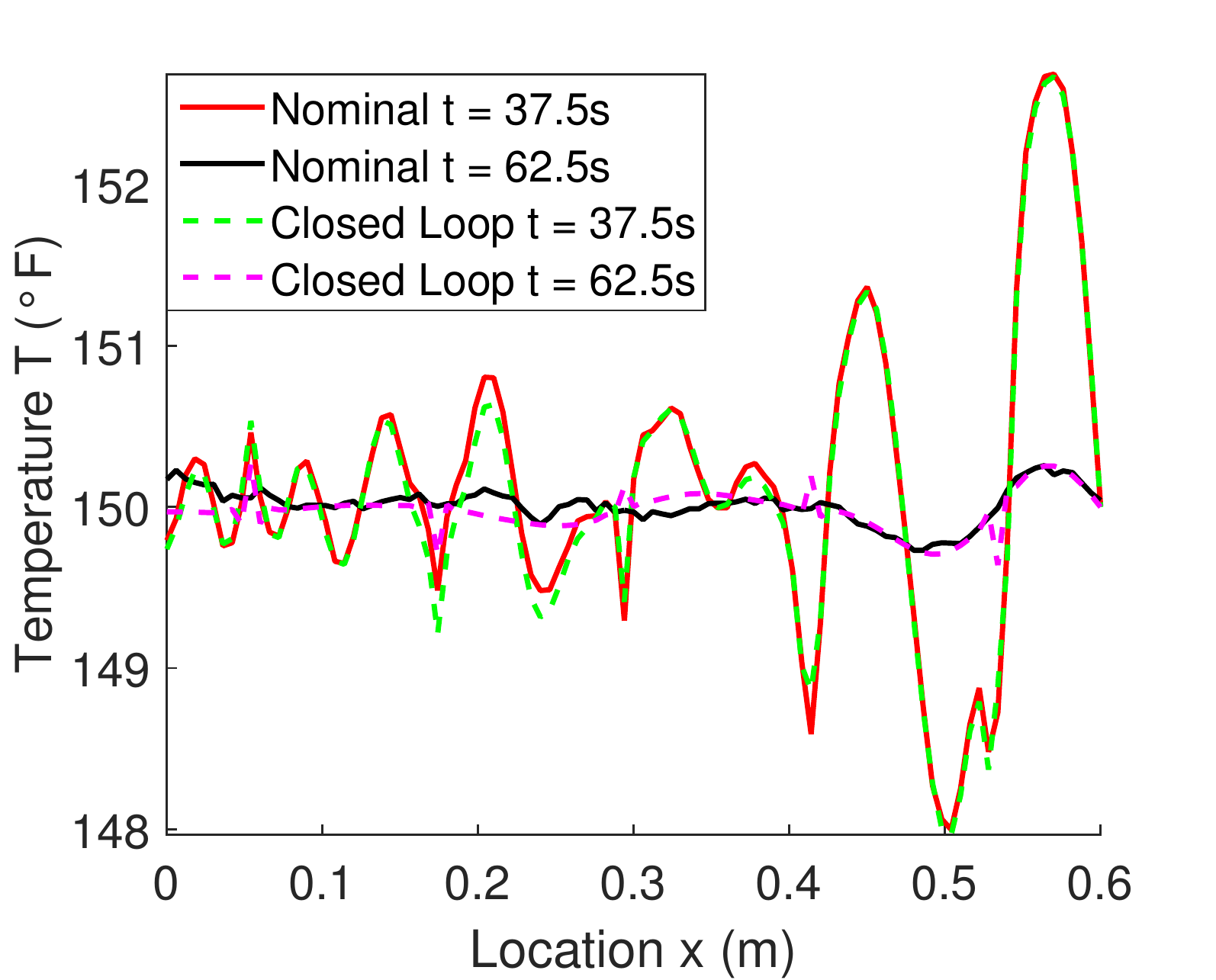}}
\subfigure[Estimation Error at  $x = 0.4L$.]{
\includegraphics[width= 0.3\textwidth]{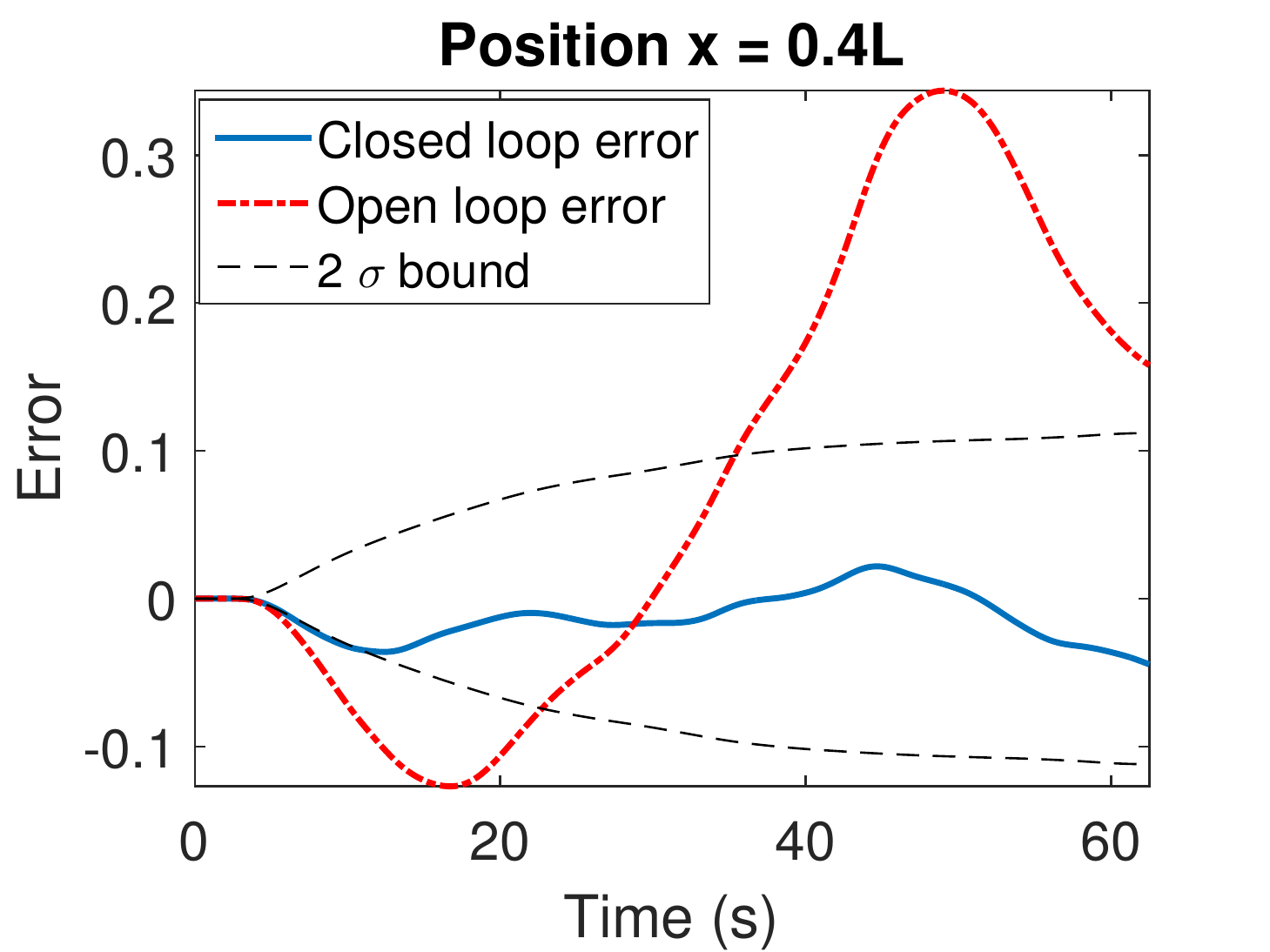}}
\subfigure[Estimation Error at $x = 0.9L$]{
\includegraphics[width= 0.3\textwidth]{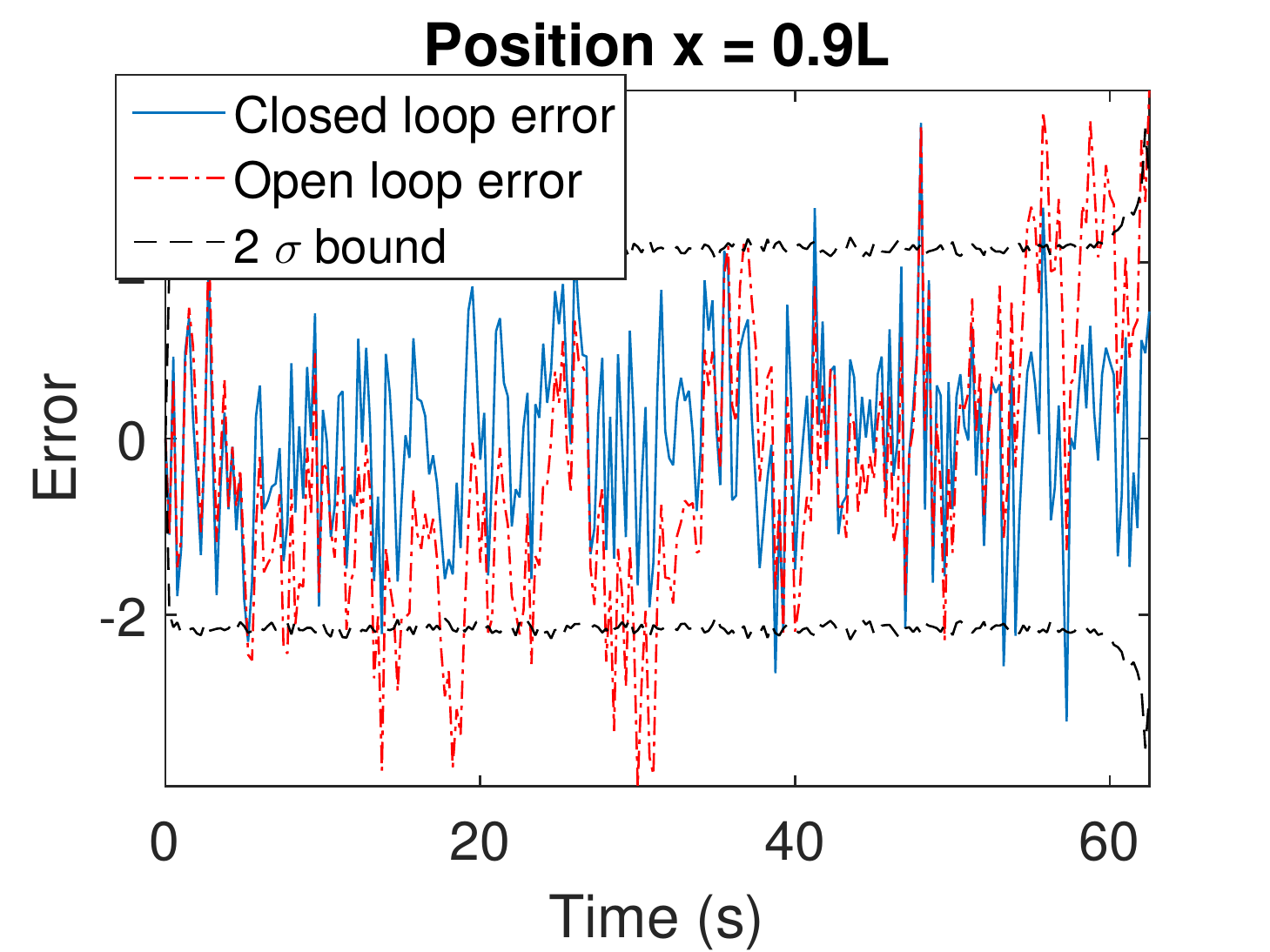}}
\caption{Performance of the Proposed Approach}\label{kf}
\end{figure*}

We solve the open loop optimization problem, and the normal (belief mean) trajectory and optimal control are shown in Fig. \ref{kf}(a) and Fig. \ref{kf}(b) respectively.

The implementation of time-varying  ERA algorithm to identify the linearized system is performed as follows. The  size of the generalized Hankel matrix $H_k^{(p, q)}$ is $pn_y \times q n_u$, and as discussed before, the design parameters $p$ and $q$ should be chosen such that $\min \{pn_y, qn_u \} \geq n_x$, which for the current problem, $n_x = 100, n_u = 5, n_y = 5$. We select $p, q$ by trial and error, i.e., we start with some initial guess of $p, q$, compare the impulse responses of the original system and the identified system, and check if the accuracy of the identified system is acceptable. Here, we choose $p = q = 15$. Therefore, the size of the generalized Hankel matrix is $75 \times 75$. The rank of the Hankel matrix is 20, and hence, the order of the identified LTV system $n_r = 20$.

We run $M$ parallel simulations to estimate the generalized Markov parameters $\{h_{k, j} \}_{j = 0}^{k}, k = 1, 2, \cdots, N -1$. We perturb the open loop optimal control $\{\bar{u}_k \}_{k = 0}^{N - 1}$ with impulse, i.e., denote $\{\delta u_k^i \}_{k = 0}^{N - 1}$ as the input perturbation sequence in the $i^{th}$ simulation, and $\{\delta u_k^i \}_{k = 0}^{N - 1} = (0, 0, \cdots, 0.01, \cdots, 0)$, where only the $i^{th}$ element is nonzero. Therefore,  we choose design parameter $M = N$. In each simulation, we collect the outputs $\{\delta y_k^i \}_{k = 0}^{N-1}$  in (\ref{noise_free_output}) corresponding to the control input $\{\bar{u}_k + \delta u_k^i \}_{k  = 0}^{N - 1}$, and solve the least squares problem using (\ref{gmarkov}).

 The rank of the Hankel matrix $n_r = 20$, and hence, the identified reduced order system  $\hat{A}_k \in \Re^{20 \times 20}$.  Due to the separation principle, the feedback design decouples into the solution of two 20 x 20 Ricatti equations, one for the controller and one for the Kalman filter: compare this to the 10100 x 10100 problem that would need to be solved if using an ILQG based approach. With the identified linearized system, we design the closed loop controller. We run 1000 individual simulations with process noise $w_k \sim N(0, W)$ and measurement noise $v_k \sim N(0, V)$,where $W = I_{5 \times 5}$, $V = I_{5 \times 5}$.

In Fig. \ref{kf}, we show the performance of the proposed approach. We calculate the identified Markov parameters using $(\hat{A}_k, \hat{B}_k, \hat{C}_k)$, and compare with the actual generalized Markov parameters (calculated using impulse responses). The Markov parameters $h_{k, j} \in \Re^{ 5 \times 5}$, and we show the comparison for one input-output channel at time step $k = 250$ in Fig.\ref{kf}(c). It can be seen that the identified LTV system using time-varying ERA approach can approximate the linearized deviation system accurately. In Fig. \ref{kf}(d), we compare the averaged closed loop trajectory with the nominal trajectory at time $t = 37.5s, t = 62.5s$. In Fig. \ref{kf}(e) - (f), we randomly choose two positions, and show  the errors between the actual trajectory and optimal trajectory with $2 \sigma$ bound in one simulation. For comparison, the open loop error is also shown in the figure. 

It can be seen that the averaged state estimates over 1000 Monte-Carlo simulations runs are close to the open loop optimal trajectory, which implies that the control objective to minimize the expected cost function could be achieved using the proposed approach. In this partially observed problem, the computationally complexity of designing the online estimator and controller using the identified ROM model are reduced by orders of $O(\frac{ n_x^4}{n_r^2}) = O(10^5)$, and for a general three dimensional problem this reduction could be even more significant.

\section{Conclusion}
In this paper, we have proposed a separation based design of the stochastic optimal control problem for systems with unknown nonlinear dynamics and partially observed states. First, we design a deterministic open-loop optimal trajectory in belief space.  Then we identify the nominal linearized system using time-varying ERA. The open-loop optimization and system identification are implemented offline, using the impulse responses of the system, and an LQG controller based on the ROM is implemented online.  The offline learning procedure is simple, and the online implementation is fast. We have tested the proposed approach on a one dimensional nonlinear heat transfer problem, and showed the performance of the proposed approach.

\appendices
\section{Ensemble Kalman Filter}\label{app_enkf}
Consider the discrete time nonlinear dynamical system (\ref{original system}) with Gaussian belief states $b_k = (\mu_k, \Sigma_k)$. Assume that $b_0$ is known.  Denote  $\{x_{k, (j)}^- \}, j = 1, \cdots, m$ as an $m$-member forecast ensemble at time step $k$, where the subscript $(j)$ denotes the $j^{th}$ member, the forecast mean $b_k^-$ and covariance $\Sigma_k^-$ are defined as
\begin{align}
&b_k^- = \frac{1}{m} \sum_{j = 1}^m x_{k, (j)}^-, \nonumber\\
&\Sigma_k^- = \frac{1}{m - 1} \sum_{i = 1}^m (x_{k, (j)}^-  - b_k^- )(x_{k, (j)}^- - b_k^-)'.
\end{align}

The measurement ensemble at time $k$ is $\{y_{k, (j)}^- \}, j = 1, \cdots, m$, where
\begin{align}
y_{k, (j)}^- = h(x_{k, (j)}^-, v_k), v_k \sim N(0, V).
\end{align}
The corresponding mean and covariance $y_k^-, P_k^y$ are defined similarly.

Denote the cross-covariance matrix $P_{k}^{xy}$ between the state and measurement ensembles at time $k$ as:
\begin{align}
 P_{k}^{xy} = \sum_{j = 1}^m (x_{k, (j)} - x_{k}^-)(y_{k, (j)} - y_{k}^-)'.
\end{align}

Given an initial belief $b_0$, and a control sequence $\{u_k \}_{k = 0}^{N - 1}$, the EnKF used to estimate the belief state $\{b_k\}_{k = 0}^N$ is summarized in Algorithm \ref{algo_enkf}.
\begin{algorithm}[!htbp]
\caption{EnKF}
 \label{algo_enkf}
  \begin{algorithmic}[1]
 \Require{Start belief $b_0$, control sequence $ \{u_k \}_{k = 0}^{N - 1}$.}
 \Ensure{Belief nominal trajectory $\{b_k \}_{k  = 0}^N$.}
 \For{$k = 0 : N$}
  \State Sample $m$ forecast ensemble members:
 \begin{align}
& x_{k + 1, (j)}^- = f(b_k, u_k, w_k),  \nonumber\\
 &y_{k+1, (j)}^- = h(x_{k+1,(j)}^-, v_{k+1}), j =  1, \cdots, m
 \end{align}
 where  $w_k \sim (0, W), v_k \sim (0, V)$, and calculate $b_{k+1}^-, \Sigma_{k+1}^- $, $y_{k+1}^-, P_{k+1}^y, P_{k+1}^{xy}$.
 \State Propagate the underlying noiseless system and take  the measurement $y_{k+1}$.
 \begin{align}
 &x_{k + 1} = f(x_k, u_k, 0), y_{k+1} = h(x_{k+1}, 0).
 \end{align}
 \State Update the posterior mean and covariance
 \begin{align}
 &b_{k + 1} = b_{k + 1}^{-} + K_e(y_{k+1} - y_{k+1}^{-}), \nonumber\\
 &\Sigma_{k + 1} = \Sigma_{k+1}^- - P_{k+1}^{xy} (P_{k+1}^{y})^{-1}P_{k+1}^{yx},
 \end{align}
 where
 \begin{align}
  K_e = P_{k+1}^{xy}(P_{k+1}^y)^{-1}.
 \end{align}
 \EndFor
\end{algorithmic}
\end{algorithm}

\section{Closed Loop Feedback Controller}\label{sec_lqg}
Given the identified linear deviation system (\ref{rom}),  the separation principle could be used, and hence, we design a Kalman filter and an LQR controller separately.

The feedback controller is:
\begin{align}
\delta u_k = -L_k \delta \hat{a}_k,
\end{align}
where $\delta \hat{ a}_k$ is the estimate from a Kalman observer, and the feedback gain $L_k$ is computed by solving two decoupled Riccati equations as follows.
\begin{align}\label{r_1}
L_k = (\hat{B}_k' S_{k + 1} \hat{B}_k + R_k)^{-1} \hat{B}_k' S_{k + 1} \hat{A}_k,
\end{align}
where $S_k$ is determined by running the following Riccati equation backward in time:
\begin{align}\label{r_2}
S_k = & \hat{A}_k' S_{k + 1} \hat{A}_k  + Q_k \nonumber\\
 &- \hat{A}_k S_{k + 1} \hat{B}_k(\hat{B}_k' S_{k + 1} \hat{B}_k + R_k)^{-1} \hat{B}_k' S_{k + 1} \hat{A}_k,
\end{align}
with terminal condition $S_N = Q_N$.

The Kalman filter observer is designed as follows:
\begin{align}\label{k_1}
\delta \hat{a}_{k + 1} = &\hat{A}_k \delta \hat{a}_k + \hat{B}_k \delta u_k  \nonumber\\
&+ K_{k + 1}(\delta y_{k + 1} - \hat{C}_{k + 1} (\hat{A}_k \delta \hat{a}_k + \hat{B}_k \delta u_k)),
\end{align}
with $\delta y_k = h(x_k, v_k) - h(\bar{\mu}_k, 0)$, and the covariance of the estimation is:
\begin{align}\label{k_2}
P_{k + 1} = &\hat{A}_k(P_k - P_k \hat{C}_k' (\hat{C}_k P_k \hat{C}_k' + V)^{-1} \hat{C}_k P_k) \hat{A}_k' \nonumber\\
&+ \hat{B}_k W \hat{B}_k',
\end{align}
where the Kalman gain is:
\begin{align}\label{k_3}
K_k = P_k \hat{C}_k' (\hat{C}_k P_k \hat{C}_k' + V)^{-1}.
\end{align}

\bibliographystyle{IEEEtran}
\bibliography{CDC_refs}

\begin{thebibliography}{10}
\providecommand{\url}[1]{#1}
\csname url@samestyle\endcsname
\providecommand{\newblock}{\relax}
\providecommand{\bibinfo}[2]{#2}
\providecommand{\BIBentrySTDinterwordspacing}{\spaceskip=0pt\relax}
\providecommand{\BIBentryALTinterwordstretchfactor}{4}
\providecommand{\BIBentryALTinterwordspacing}{\spaceskip=\fontdimen2\font plus
\BIBentryALTinterwordstretchfactor\fontdimen3\font minus
  \fontdimen4\font\relax}
\providecommand{\BIBforeignlanguage}[2]{{%
\expandafter\ifx\csname l@#1\endcsname\relax
\typeout{** WARNING: IEEEtran.bst: No hyphenation pattern has been}%
\typeout{** loaded for the language `#1'. Using the pattern for}%
\typeout{** the default language instead.}%
\else
\language=\csname l@#1\endcsname
\fi
#2}}
\providecommand{\BIBdecl}{\relax}
\BIBdecl

\bibitem{dp_bertsekas}
D.~P. Bertsekas, \emph{Dynamic Programming and Optimal Control, Two Volume
  Set}, 2nd~ed.\hskip 1em plus 0.5em minus 0.4em\relax Athena Scientific, 1995.

\bibitem{dp_num}
M.~Falcone, ``{Recent Results in the Approximation of Nonlinear Optimal Control
  Problems},'' in \emph{Large-Scale Scientific Computing LSSC}, 2013.

\bibitem{bryson}
A.~Bryson and H.~Y.-C., \emph{Applied Optimal Control: Optimization, Estimation
  and Control}.\hskip 1em plus 0.5em minus 0.4em\relax Washington: Hemisphere
  Pub. Corp., 1975.

\bibitem{ddp}
D.~Jacobsen and D.~Mayne, \emph{Differential Dynamic Programming}.\hskip 1em
  plus 0.5em minus 0.4em\relax Elsevier, 1970.

\bibitem{sddp}
E.~Theoddorou, Y.~Tassa, and E.~Todorov, ``{Stochastic Differential Dynamic
  Programming},'' in \emph{Proceedings of American Control Conference}, 2010.

\bibitem{ilqg1}
E.~Todorov and W.~Li, ``{A generalized iterative LQG method for locally-optimal
  feedback control of constrained nonlinear stochastic systems},'' in
  \emph{Proceedings of American Control Conference}, 2005, pp. 300 -- 306.

\bibitem{ilqg2}
W.~Li and E.~Todorov, ``Iterative linearization methods for approximately
  optimal control and estimation of non-linear stochastic system,''
  \emph{International Journal of Control}, vol.~80, no.~9, pp. 1439--1453,
  2007.

\bibitem{adp}
R.~P. Bithmead, V.~Wertz, and M.~Gerers, \emph{Adaptive Optimal Control: The
  Thinking Man's G.P.C}.\hskip 1em plus 0.5em minus 0.4em\relax Prentice Hall
  Professional Technical Reference, 1991.

\bibitem{adp_3}
X.~Zhong, H.~He, H.~Zhang, and Z.~Wang, ``{Optimal Control for Unknown
  Diiscrete-Time Nonlinear Markov Jump Systems Using Adaptive Dynamic
  Programming},'' \emph{IEEE Transactions on Neural networks and learning
  systems}, vol.~25, no.~12, pp. 2141--2155, 2014.

\bibitem{rl_1}
S.~G. Khan \emph{et~al.}, ``{Reinforcement learning and optimal adaptive
  control: An overview and implementation examples},'' \emph{Annual Reviews in
  Control}, vol.~36, pp. 42--59, 2012.

\bibitem{rl_3}
D.~Mitrovic, S.~Klanke, and S.~Vijayakumar, \emph{Adaptive Optimal Feedback
  Control with Learned Internal Dynamics Models}.\hskip 1em plus 0.5em minus
  0.4em\relax Berlin, Heidelberg: Springer Berlin Heidelberg, 2010, pp. 65--84.

\bibitem{RLHD4}
S.~Levine and P.~Abbeel, ``{Learning Neural Network Policies with Guided Search
  under Unknown Dynamics},'' in \emph{Advances in Neural Information Processing
  Systems}, 2014.

\bibitem{RLHD5}
S.~Levine and K.~Vladlen, ``{Learning Complex Neural Network Policies with
  Trajectory Optimization},'' in \emph{Proceedings of the International
  Conference on Machine Learning}, 2014.

\bibitem{RLHD1}
R.~Akrour, A.~Abdolmaleki, H.~Abdulsamad, and G.~Neumann, ``{Model Free
  Trajectory Optimization for Reinforcement Learning},'' in \emph{Proceedings
  of the International Conference on Machine Learning}, 2016.

\bibitem{RLHD2}
E.~Todorov and Y.~Tassa, ``{Iterative Local Dynamic Programming},'' in
  \emph{Proc. of the IEEE Int. Symposium on ADP and RL.}, 2009.

\bibitem{baxter}
J.~Baxter and P.~Bartlett, ``{Infinite Horizon Policy-Gradient Estimation},''
  \emph{Journal of Artificial Intelligence Research}, vol.~15, pp. 319--350,
  2001.

\bibitem{sutton3}
R.~S. Sutton, D.~Mcallester, S.~Singh, and Y.~Mansour, ``{Policy Gradient
  Methods for Reinforcement Learning with Function Approximation},'' in
  \emph{Proc. 1999 Neural Information Proc. Sys.}, 1999.

\bibitem{marbach}
P.~Marbach, \emph{{Simulation based Optimization of Markov Reward Processes,
  PhD Thesis}}.\hskip 1em plus 0.5em minus 0.4em\relax Boston, MA:
  Massachusetts Institute of Technology, 1999.

\bibitem{Survey_PS}
M.~P. Deisenroth, G.~Neumann, and J.~Peters, ``{A Survey on Policy Search for
  Robotics},'' in \emph{Foundations and Trends in Robotics}, 2013, pp. 1--142.

\bibitem{tv_era}
M.~Majji, J.-N. Juang, and J.~L. Junkins, ``{Time-varying Eigensystem
  Realization Algorithm},'' \emph{Journal of Guidance, Control, and Dynamics},
  vol.~33, no.~1, pp. 13--28, 2010.

\bibitem{enkf1}
S.Gillijins \emph{et~al.}, ``{What Is the Ensemble Kalman Filter and How Well
  Does it Work?}'' in \emph{Proceedings of the 2006 American Control
  Conference}, 2006, pp. 4448--4453.

\bibitem{Separation}
M.~Rafieisakhaei, S.~Chakravorty, and P.~R. Kumar, ``{A Near-Optimal Separation
  Principle for Nonlinear Stochastic Systems Arising in Robotic Path Planning
  and Control},'' in \emph{56$^{th}$ IEEE Conference on Decision and
  Control(CDC)}, 2017.

\bibitem{platt10}
R.~Platt, R.~Tedrake, L.~Kaelbling, and T.~Lozano-Perez, ``Belief space
  planning assuming maximum likelihood observatoins,'' in \emph{Proceedings of
  Robotics: Science and Systems (RSS)}, June 2010.

\bibitem{enkf2}
F.~Hamilton, T.~Berry, and T.~Sauer, ``{Ensemble Kalman Filtering without a
  Model},'' \emph{Physical Review X}, vol.~6, p. 011021, 2016.

\bibitem{gradient}
A.E.Bryson and W.~Denham, ``A steepest-ascent method for solving optimum
  programming problems,'' \emph{Journal of Applied Mechanics}, vol.~29, no.~2,
  1962.

\bibitem{sim_opt}
A.~Gosavi, \emph{Simulation-based optimization: Parametric optimization
  techniques and reinforcement learning}.\hskip 1em plus 0.5em minus
  0.4em\relax Norwell, MA, USA: Kluwer Academic Publishers, 2003.

\bibitem{antoulas}
A.~Antoulas, \emph{{Approximation of Large Scale Dynamical Systems}}.\hskip 1em
  plus 0.5em minus 0.4em\relax Philadelphia, PA: SIAM, 2005.

\bibitem{skelton1}
A.M.King, U.B.Desai, and R.E.Skelton, ``{A Generalized Approach to q-Markov
  Covariance Equivalent Realizations for Discrete Systems},''
  \emph{Automatica}, vol.~24, no.~4, pp. 507--515, 1988.

\bibitem{dlqg_skelton}
G.~Shi and R.~E. Skelton, ``{Markov Data-Based LQG Control},'' \emph{Journal of
  Dynamic Systems, Measurement, and Control}, vol. 122, pp. 551--559, 2000.

\bibitem{dlqg_2}
W.~Favoreel \emph{et~al.}, ``Closed-loop model-free subspace-based
  lqg-design,'' in \emph{Proccedings of the 7th Mediterranean Conference on
  Control and Automation}, 1999, pp. 1926--1939.

\bibitem{dd_1}
Z.-S. Hou and Z.~Wang, ``From model-based control to data-driven control:
  Survey, classification and perspective,'' \emph{Information Sciences}, vol.
  235, no. 3-35, 2013.

\bibitem{dd_2}
X.~Wu, J.~Shen, Y.~Li, and K.~Y.Lee, ``{Data-Driven Modeling and Predictive
  Control for Bioler-Turbine Unit},'' \emph{IEEE Transactions on energy
  conversion}, vol. 228, no.~3, pp. 470--481, 2013.

\end{thebibliography}

\end{document}